\documentclass[11pt]{amsart}

\usepackage[a4paper,margin=28mm]{geometry}
\usepackage{amsmath,amssymb,amsfonts,amsthm,mathrsfs}
\usepackage{mathtools}
\usepackage{microtype}
\usepackage{hyperref}
\usepackage{xcolor}
\hypersetup{
    colorlinks,
    linkcolor={red!50!black},
    citecolor={blue!100!black},
    urlcolor={blue!100!black}
}

\numberwithin{equation}{section}

\newtheorem{theorem}{Theorem}[section]
\newtheorem{lemma}[theorem]{Lemma}

\newcommand{\e}{\mathrm{e}}
\newcommand{\T}{\mathbb{T}}
\newcommand{\dd}{\mathrm{d}}
\newcommand{\1}{\mathbf{1}}
\newcommand{\norm}[1]{\left\lVert #1\right\rVert}

\title[Liouville exponential sums]
{On the lower bound for the $L^1$ norm of Liouville exponential sums. I.}

\author[J. Li]{Jiangpeng Li}
\author[Y.-C. Sun]{Yu-Chen Sun}
\author[H. Wang]{Hui Wang}

\address
{School of Mathematics, Shandong University, Jinan 250100, China}
\email{jpli@mail.sdu.edu.cn}
\address
{Department of Mathematics, University of Bristol, Woodland Rd, Bristol BS8 1UG}
\email{yuchensun93@163.com}
\address
{The University of Hong Kong Shenzhen Institute of Research and Innovation,
Shenzhen, Guangdong 518057, China}
\email{wh0315@mail.sdu.edu.cn}

\makeatletter
\@ifundefined{subjclassname@2020}{%
  \@namedef{subjclassname@2020}{\textup{2020} Mathematics Subject Classification}}{}
\makeatother
\date{}

\begin{document}

\begin{abstract}
Let $\lambda$ be the Liouville function. We prove
\[
\mathcal{I}_{\lambda}(X)=\int_{\T}\left|\sum_{1\leqslant n\leqslant X}\lambda\left(n\right)\e\left(n\alpha\right)\right|
 \dd\alpha
 \gg X^{\frac{1}{4}}.
\]
This removes the $X^\varepsilon$ loss in the
previous result of Pandey and Radziwi\l\l~\cite{PandeyRadziwill2026}. In the second paper in this series, we will prove that
\[
\limsup_{X \to \infty} \frac{\log \mathcal{I}_{\lambda}(X)} {\log X} \ge \frac{1}{3}.
\]
\end{abstract}

\maketitle

\section{Introduction}

The $L^1$ norms of exponential sums are a classical subject in harmonic
analysis. We write $\e\left(u\right)=\exp\left(2\pi\mathrm{i}u\right)$ and
$\T=\mathbb R/\mathbb Z$, equipped with Haar measure of total mass one. For a
complex sequence $a\left(n\right)$ and a real number $N\geqslant1$, define
\begin{equation}\label{eq:I-definition}
 \mathcal I_a\left(N\right)
 =\int_{\T}\left|\sum_{1\leqslant n\leqslant N}a\left(n\right)\e\left(n\alpha\right)\right|\dd\alpha.
\end{equation}
Konyagin \cite{Konyagin1982} and McGehee, Pigno and Smith
\cite{McGeheePignoSmith1981} established Littlewood's logarithmic lower bound
for exponential sums with coefficients of modulus at least one. Stronger lower bounds are available for certain arithmetic coefficients. For example, let $\mu_k$ be the characteristic function of the $k$-free integers, where $k\geqslant2$ is fixed. Balog and Ruzsa
\cite{BalogRuzsa2001} proved
$\mathcal I_{\mu_k}\left(X\right)\asymp_k X^{\frac{1}{k+1}}$.
In particular, $\mu_2=\mu^2$, where $\mu$ denotes the M\"obius function. By \eqref{convolution}, this implies
\(\mathcal I_\mu\left(X\right)\gg X^{\frac{1}{6}}\) \cite{BalogRuzsa2001}. Their approach does not seem to apply directly to the Liouville function, defined by
$\lambda\left(n\right)=\left(-1\right)^{\Omega\left(n\right)}$. Here $\Omega\left(n\right)$ denotes the number of prime factors of $n$, counted with multiplicity.

For Liouville coefficients, Littlewood's inequality gives the lower bound,
while the Cauchy--Schwarz inequality and Parseval's identity give the upper
bound
\begin{equation}\label{eq:general-bounds}
 \log\left(N+1\right)\ll\mathcal I_\lambda\left(N\right)\ll\sqrt N.
\end{equation}
Balog and Perelli \cite{BalogPerelli1998} obtained
$\mathcal I_\lambda\left(X\right)\gg X^{\frac{c}{\log\log X}}$
for sufficiently large $X$ and an absolute constant $c>0$. More recently, Pandey and Radziwi\l\l~
\cite{PandeyRadziwill2026} proved that, for every
$\varepsilon>0$,
\begin{equation}\label{eq:PR-prefix}
 \mathcal I_{\lambda}\left(X\right)\gg_\varepsilon X^{\frac{1}{4}-\varepsilon}.
\end{equation}
They also observed that the factor $X^{-\varepsilon}$ in
\eqref{eq:PR-prefix} can be replaced by a negative power of $\log X$.
We remove this loss.
\begin{theorem}\label{thm:prefix}
Let $X\geqslant2$. Then
\begin{equation}\label{eq:prefix}
 \mathcal I_\lambda\left(X\right)
 =\int_{\T}\left|\sum_{1\leqslant n\leqslant X}\lambda\left(n\right)\e\left(n\alpha\right)\right|
 \dd\alpha
 \gg X^{\frac{1}{4}}.
\end{equation}
\end{theorem}

The proof of Theorem~\ref{thm:prefix} is essentially that of Pandey and
Radziwi\l\l~\cite{PandeyRadziwill2026}. They introduce a suitable smooth
weight and start from the convolution identity \eqref{convolution}; see
\cite[Lemma~3.2]{PandeyRadziwill2026}. They then bound the $L^1$ norm
from below by integrating near rational points with prime denominators
up to $\sqrt X$ and summing over the reduced numerators to obtain
Ramanujan sums. The key point of their argument is the evaluation of
the resulting smoothed Liouville sums in
\cite[Proposition~2.2]{PandeyRadziwill2026}, using a kernel whose Mellin
transform contains the relevant $L$-function.

We use this argument with two modifications. In
Lemma~\ref{lem:weight}, we apply their kernel construction, going back
to Pintz~\cite{Pintz1982}, to the pole of $\zeta(2s)$ at
$s=\frac{1}{2}$ in place of a zero of a Dirichlet $L$-function.
This gives a fixed real-valued function
$w\in C_c^\infty(\mathbb R)$, supported in $[0,\frac{1}{2}]$, such that
\[
 F(Y):=\sum_{n\geqslant1}\lambda(n)w(n/Y)
 =Y^{\frac{1}{2}}+O_A(Y^{-A})
 \qquad(A>0,\ Y\geqslant1).
\]
We give a direct proof of this form of the construction using the identity
\begin{equation}\label{eq:lambda-square}
 \sum_{d\mid r}\lambda(d)=\1_{\{r\text{ is a square}\}}
 \qquad(r\geqslant1)
\end{equation}
and Poisson summation.

We also extend the sum over prime denominators to all denominators
$q\leqslant Q$, where $Q\asymp\sqrt X$.
After multiplication by $\lambda(q)$, each Ramanujan-weighted Liouville
sum is a linear combination of the values $F(X/d)$, $d\mid q$, with
nonnegative coefficients; see \eqref{eq:positive-combination}.
For sufficiently large $X$, these values are positive uniformly for
$q\leqslant Q$, and the term $d=q$ gives a lower bound
$\gg\sqrt{Xq}$. Summing over all such $q$ and applying the $L^1$
large sieve inequality in Lemma~\ref{lem:sampling} yields
Theorem~\ref{thm:prefix} without $X^{\varepsilon}$ loss.

The construction in Appendix~\ref{app:weight} is essentially that of
Pandey and Radziwi\l\l~\cite[Proposition~2.2 and Lemma~4.1]{PandeyRadziwill2026}.
For any fixed zero $\rho=\beta+\mathrm{i}\gamma$ of $\zeta$ with
$\frac{1}{2}\leqslant\beta<1$ and $\gamma\ne0$, combining the weight in
Lemma~\ref{lem:zero-weight} with our sampling argument over all
denominators gives
\[
 \mathcal I_\lambda(X)\gg_\rho X^{\frac{\beta}{2}}
 \qquad(X\geqslant2),
\]
with no $X^\varepsilon$ loss.

\subsection*{Acknowledgements}

The authors thank Maksym Radziwi\l\l\ for helpful comments.

\subsection*{AI assistance}

The authors used ChatGPT (OpenAI) to assist with the proof of Lemma~\ref{lem:weight} and the argument in Appendix~\ref{app:weight}, the search for references and the preparation of the manuscript. The authors take full responsibility for the results and their presentation.

\subsection*{Notation}

For a trigonometric polynomial $P$ on $\T$, write
$\norm{P}_{L^1\left(\T\right)}=\int_{\T}|P\left(\alpha\right)|\dd\alpha$.
For $f\in L^1\left(\mathbb R\right)$, we use the Fourier transform convention
$\widehat f\left(\xi\right)=\int_{\mathbb R}f\left(u\right)\e\left(-u\xi\right)\dd u$
for $\xi\in\mathbb R$.
For a trigonometric polynomial $P$, its Fourier coefficients are
$\widehat P\left(n\right)=\int_{\T}P\left(\alpha\right)\e\left(-n\alpha\right)\dd\alpha$ for $n\in\mathbb Z$.
For trigonometric polynomials $P$ and $Q$, convolution on $\T$ is taken with
respect to Haar measure and is given by
\begin{equation} \label{convolution}\left(P*Q\right)\left(\alpha\right)=\int_{\T}P\left(\alpha-\beta\right)Q\left(\beta\right)\dd\beta.
\end{equation}
The symbol $\1_E$ denotes the indicator of a set $E$, and
$C_c^\infty\left(U\right)$ denotes the space of smooth functions compactly supported in an open
set $U\subseteq\mathbb R$.
\section{Two auxiliary lemmas}

Lemma~\ref{lem:sampling} is the $p=1$ case of the $L^p$ large sieve
inequality of Golinskii, Lubinsky and Nevai~\cite[Corollary~1.2]{GolinskiiLubinskyNevai2001}.
Lemma~\ref{lem:weight} applies the kernel construction of Pandey and
Radziwi\l\l~\cite[Proposition~2.2]{PandeyRadziwill2026}, going back to
Pintz~\cite{Pintz1982}, to the pole of $\zeta(2s)$ at $s=\frac{1}{2}$.
The proof below uses the square-detection identity of
Yakubovich~\cite[Theorem~5]{Yakubovich2015} and Poisson summation.
We include full proofs for completeness.

\begin{lemma}\label{lem:sampling}
Let $N\geqslant1$ be an integer and $0<\delta\leqslant1$. Let
$\mathcal R\subset\T$ be a finite set whose distinct elements have circular
distance at least $\delta$. If
$T\left(\alpha\right)=\sum_{|n|\leqslant N}b_n\e\left(n\alpha\right)$, then
\begin{equation}\label{eq:sampling}
 \sum_{r\in\mathcal R}|T\left(r\right)|
 \ll \left(N+\delta^{-1}\right)\int_{\T}|T\left(\alpha\right)|\dd\alpha.
\end{equation}
\end{lemma}

\begin{proof}
For each integer $M\geqslant1$, let $F_M$ denote the Fej\'er kernel
\begin{align*}
 F_M\left(t\right)=\frac{1}{M}\left|\sum_{j=0}^{M-1}\e\left(jt\right)\right|^2,
 \qquad K_N=2F_{2N}-F_N.
\end{align*}
Since $\widehat F_M\left(n\right)=\max\{0,1-\frac{|n|}{M}\}$, we have
$\widehat K_N\left(n\right)=1$ for $|n|\leqslant N$, and hence $T=T*K_N$. By the finite
geometric-series formula, we have
\begin{equation}\label{eq:kernel}
 |K_N\left(t\right)|\ll\frac{N}{\left(1+N\norm{t}\right)^2},
\end{equation}
where $\norm{t}$ is the distance from $t$ to the nearest integer.

Fix $\beta\in\T$. For each integer $j\geqslant0$, the points $r\in\mathcal R$
satisfying $\frac{j}{N}\leqslant\norm{r-\beta}<\frac{j+1}{N}$ lie in the union of at most two arcs of length at most $\frac{1}{N}$. Since distinct
points have circular distance at least $\delta$, there are
$O\left(1+\left(N\delta\right)^{-1}\right)$ such points. Hence
\begin{align*}
 \sum_{r\in\mathcal R}|K_N\left(r-\beta\right)|
 \ll N\left(1+\frac{1}{N\delta}\right)
       \sum_{j\geqslant0}\left(1+j\right)^{-2}
 \ll N+\delta^{-1}.
\end{align*}
Using $T=T*K_N$, the triangle inequality and the last estimate, we obtain
\eqref{eq:sampling}.
\end{proof}

We next construct a fixed smooth weight for which the Liouville sum has an explicit main term.

\begin{lemma}\label{lem:weight}
There is a fixed real-valued function $w\in C_c^\infty\left(\mathbb R\right)$, supported
in $[0,\frac{1}{2}]$, such that, for every $A>0$ and $Y\geqslant1$, we have
\begin{equation}\label{eq:weighted-sum}
 F\left(Y\right):=\sum_{n\geqslant1}\lambda\left(n\right)w\left(\frac{n}{Y}\right)
 =Y^{\frac{1}{2}}+O_A\left(Y^{-A}\right).
\end{equation}
\end{lemma}

\begin{proof}
Choose a fixed nonzero nonnegative function
$\eta\in C_c^\infty\left(\left(\frac{1}{4},\frac{1}{2}\right)\right)$ and put
\begin{align*}
 V\left(u\right)=4\left(\int_0^\infty\eta\left(v\right)v^{-\frac{3}{2}}\dd v\right)^{-1}
 \frac{\dd\eta\left(u\right)}{\dd u}.
\end{align*}
Integrating by parts, we obtain
\begin{equation}\label{eq:V-moments}
 \int_0^\infty V\left(u\right)\dd u=0,
 \qquad
 \frac{1}{2}\int_0^\infty V\left(u\right)u^{-\frac{1}{2}}\dd u=1.
\end{equation}
For $u>0$, put $w\left(u\right)=\sum_{m\geqslant1}V\left(mu\right)$, and set $w\left(u\right)=0$ for
$u\leqslant0$. The sum is locally finite on $\left(0,\infty\right)$, and $w\left(u\right)=0$ for
$u\geqslant\frac{1}{2}$.

Fix an integer $j\geqslant0$ and put $U_j\left(t\right)=t^jV^{\left(j\right)}\left(t\right)$, extended by
zero to the real line. By repeated integration by parts and
\eqref{eq:V-moments}, we have
\begin{align*}
 \int_{\mathbb R}U_j\left(t\right)\dd t
 =\left(-1\right)^j j!\int_{\mathbb R}V\left(t\right)\dd t=0.
\end{align*}
For $u>0$, by the support of $U_j$, we have
$u^jw^{\left(j\right)}\left(u\right)=\sum_{m\in\mathbb Z}U_j\left(mu\right)$. The integral identity above gives $\widehat U_j\left(0\right)=0$. For $0<u\leqslant1$,
applying the Poisson summation formula, we obtain
\begin{align*}
 u^jw^{\left(j\right)}\left(u\right)=\frac{1}{u}\sum_{h\ne0}\widehat U_j\left(\frac{h}{u}\right).
\end{align*}
Since $U_j\in C_c^\infty\left(\mathbb R\right)$, the rapid decay of
$\widehat U_j$ gives, for every $B>1$,
\begin{align*}
 |w^{\left(j\right)}\left(u\right)|
 \ll_{B,j}u^{-j-1}\sum_{h\ne0}\left(\frac{|h|}{u}\right)^{-B}
 \ll_{B,j}u^{B-j-1}.
\end{align*}
Since $B$ can be arbitrarily large, every derivative of $w$ tends to zero as $u\to0^+$. Hence the zero extension of $w$ belongs to $C_c^\infty\left(\mathbb R\right)$.

It remains to prove \eqref{eq:weighted-sum}. Expanding the definition
of $w$ and using \eqref{eq:lambda-square}, we obtain
\begin{equation}\label{eq:square-detection}
 F\left(Y\right)=\sum_{m,n\geqslant1}\lambda\left(n\right)V\left(\frac{mn}{Y}\right)
 =\sum_{k\geqslant1}V\left(\frac{k}{Y}\right)\sum_{n\mid k}\lambda\left(n\right)
 =\sum_{\ell\geqslant1}V\left(\frac{\ell^2}{Y}\right).
\end{equation}
Put $f\left(t\right)=V\left(t^2\right)$. Then $f\in C_c^\infty\left(\mathbb R\right)$ is even and $f\left(0\right)=0$.
By \eqref{eq:V-moments}, we have
\begin{align*}
 \widehat f\left(0\right)=\int_{\mathbb R}V\left(t^2\right)\dd t
 =\int_0^\infty V\left(u\right)u^{-\frac{1}{2}}\dd u=2.
\end{align*}
Since $f$ is even and $f\left(0\right)=0$, \eqref{eq:square-detection} implies
that $F\left(Y\right)=\frac{1}{2}\sum_{\ell\in\mathbb Z}f\left(\frac{\ell}{\sqrt Y}\right)$. By the Poisson
summation formula, we have
\begin{align*}
 F\left(Y\right)=\frac{\sqrt Y}{2}\sum_{h\in\mathbb Z}\widehat f\left(h\sqrt Y\right)
 =\sqrt Y+O_A\left(Y^{-A}\right).
\end{align*}
Choose an integer $B>2A+1$. By the rapid decay of $\widehat f$, we have
\begin{align*}
 \sqrt Y\sum_{h\ne0}\left|\widehat f\left(h\sqrt Y\right)\right|
 \ll_B Y^{\frac{1-B}{2}}\sum_{h\ne0}|h|^{-B}
 \ll_A Y^{-A}.
\end{align*}
This proves \eqref{eq:weighted-sum}.
\end{proof}

\section{Proof of Theorem~\ref{thm:prefix}}

Let $w$ be as in Lemma~\ref{lem:weight}, and put
\begin{align*}
 S_X\left(\alpha\right)=\sum_{1\leqslant n\leqslant X}\lambda\left(n\right)\e\left(n\alpha\right),
 \qquad
 T_X\left(\alpha\right)=\sum_{n\geqslant1}\lambda\left(n\right)w\left(\frac{n}{X}\right)\e\left(n\alpha\right).
\end{align*}
Put $P_X\left(\alpha\right)=\sum_{n\in\mathbb Z}w\left(\frac{n}{X}\right)\e\left(n\alpha\right)$. Since $w$ is
supported in $[0,\frac{1}{2}]$ and $w\left(0\right)=0$, we have $T_X=S_X*P_X$. By \eqref{convolution}, we obtain
\[
\norm{T_X}_{L^1\left(\T\right)} \leq I_\lambda\left(X\right) \cdot \norm{P_X}_{L^1\left(\T\right)}.
\]

Poisson summation and the triangle inequality give
\begin{align*}
 P_X\left(\alpha\right)=X\sum_{k\in\mathbb Z}\widehat w\left(X\left(k-\alpha\right)\right),
 \qquad \norm{P_X}_{L^1\left(\T\right)}\leqslant\norm{\widehat w}_{L^1\left(\mathbb R\right)}.
\end{align*}
Combining these estimates gives
\begin{equation}\label{eq:multiplier}
 \norm{T_X}_{L^1\left(\T\right)}\ll \mathcal I_\lambda\left(X\right).
\end{equation}
It remains to bound $\norm{T_X}_{L^1\left(\T\right)}$ from below. Let $Q \asymp \sqrt{X
}$ and
\[
 \mathcal R_Q=
 \left\{\frac{a}{q}\bmod 1:1\leqslant q\leqslant Q,\
 0\leqslant a<q,\ \left(a,q\right)=1\right\}.
\]

Distinct elements of $\mathcal R_Q$ have circular distance at least
$Q^{-2}$ on $\T$. Applying Lemma~\ref{lem:sampling} with $N=\lceil X\rceil$
and $\delta=Q^{-2}$, we have
\[
\sum_{r\in\mathcal R_Q}|T_X\left(r\right)|
 \ll \left(\lceil X\rceil+Q^2\right)\norm{T_X}_{L^1\left(\T\right)}
 \ll X\mathcal I_\lambda\left(X\right).
\]

For each positive integer $q$, define
\begin{align*}
 c_q\left(n\right)=\sum_{\substack{0\leqslant a<q\\\left(a,q\right)=1}}\e\left(\frac{an}{q}\right),
 \qquad
 F_q\left(X\right)=\sum_{n\geqslant1}\lambda\left(n\right)c_q\left(n\right)w\left(\frac{n}{X}\right).
\end{align*}
By definition,
\begin{align*}
 F_q\left(X\right)=\sum_{\substack{0\leqslant a<q\\\left(a,q\right)=1}}T_X\left(\frac{a}{q}\right).
\end{align*}
The Ramanujan sum identity
\cite[Theorem~4.1, (4.7)]{MontgomeryVaughan2007} gives
\begin{align*}
 c_q\left(n\right)=\sum_{d\mid\left(q,n\right)}d\mu\left(\frac{q}{d}\right).
\end{align*}
Substituting this identity into the definition of $F_q\left(X\right)$, interchanging
the finite sums and using the complete multiplicativity of $\lambda$,
we obtain
\begin{align*}
 F_q\left(X\right)=\sum_{d\mid q}d\mu\left(\frac{q}{d}\right)\lambda\left(d\right)F\left(\frac{X}{d}\right).
\end{align*}
For $d\mid q$, since $\lambda\left(d\right)^2=1$, we have
$\lambda\left(q\right)\lambda\left(d\right)=\lambda\left(\frac{q}{d}\right)$. Also,
$\mu\left(\frac{q}{d}\right)\lambda\left(\frac{q}{d}\right)=\mu\left(\frac{q}{d}\right)^2$. Hence
\begin{equation}\label{eq:positive-combination}
 \lambda\left(q\right)F_q\left(X\right)
 =\sum_{d\mid q}d\,\mu\left(\frac{q}{d}\right)^2F\left(\frac{X}{d}\right).
\end{equation}

By Lemma~\ref{lem:weight}, every term on the
right side of \eqref{eq:positive-combination} is nonnegative. Thus
$\lambda\left(q\right)F_q\left(X\right)\geqslant0$. Keeping only the term $d=q$, we obtain
\begin{equation}\label{eq:Fq-lower}
 |F_q\left(X\right)|=\lambda\left(q\right)F_q\left(X\right)
 \geqslant qF\left(\frac{X}{q}\right)\geqslant\frac{1}{2}\sqrt{Xq}.
\end{equation}
Since $Q\asymp\sqrt X$, we have
\begin{equation}\label{eq:arithmetic-lower}
 \sum_{q\leqslant Q}|F_q\left(X\right)|
 \gg \sqrt X\sum_{q\leqslant Q}\sqrt q
 \gg X^{\frac{5}{4}}.
\end{equation}
 By the triangle inequality and \eqref{eq:multiplier},
we obtain
\begin{align*}
 X^{\frac{5}{4}} \ll \sum_{q\leqslant Q}|F_q\left(X\right)|
\leqslant\sum_{r\in\mathcal R_Q}|T_X\left(r\right)|
 \ll X\mathcal I_\lambda\left(X\right),
\end{align*}
which yields Theorem \ref{thm:prefix}.
\qed

\appendix
\section{A smooth weight associated with a zeta zero}\label{app:weight}

Mellin kernels that isolate a prescribed zero appear in
Pintz~\cite[Corollary~2.3, Sections~3--4]{Pintz2022} and
Pandey and Radziwi\l\l~\cite[Proposition~2.2, Lemma~4.1, Section~4]{PandeyRadziwill2026}.
Here we construct a fixed smooth compactly supported weight $w_\rho$
satisfying \eqref{eq:zero-weighted-sum} for every $A>0$.
We give the full proof for completeness.

For a function $f$ on $\left(0,\infty\right)$, write
$\mathcal M f\left(s\right)=\int_0^\infty f\left(u\right)u^{s-1}\dd u$
whenever the integral converges.

\begin{lemma}\label{lem:zero-weight}
Fix a zero $\rho=\beta+\mathrm{i}\gamma$ of $\zeta$ with
$\frac{1}{2}\leqslant\beta<1$ and $\gamma\ne0$. There is a fixed
complex-valued function $w_\rho\in C_c^\infty\left(\mathbb R\right)$,
supported in $[0,\frac{1}{2}]$, such that, for every $A>0$ and
$Y\geqslant1$,
\begin{equation}\label{eq:zero-weighted-sum}
 \sum_{n\geqslant1}\lambda\left(n\right)w_\rho\left(\frac{n}{Y}\right)
 =\zeta\left(2\rho\right)Y^\rho+O_{A,\rho}\left(Y^{-A}\right).
\end{equation}
Moreover, $\zeta\left(2\rho\right)\ne0$.
\end{lemma}

\begin{proof}
Choose a nonzero nonnegative function
$\eta\in C_c^\infty\left(\left(\frac{1}{4},\frac{1}{2}\right)\right)$
such that $\mathcal M\eta\left(\rho\right)\ne0$. For this, take the support in an interval on which $\gamma\log u$ varies by
less than $\frac{\pi}{2}$. Put $D=u\frac{\dd}{\dd u}$ and define
\[
 V=\frac{\left(D+1\right)\left(D+\frac{1}{2}\right)\eta}
 {\left(1-\rho\right)\left(\frac{1}{2}-\rho\right)
  \mathcal M\eta\left(\rho\right)}.
\]
All auxiliary functions are fixed once $\rho$ is chosen, independently
of $A$ and $Y$. Integration by parts gives
\[
 \mathcal M V\left(s\right)
 =\frac{\left(1-s\right)\left(\frac{1}{2}-s\right)
        \mathcal M\eta\left(s\right)}
 {\left(1-\rho\right)\left(\frac{1}{2}-\rho\right)
  \mathcal M\eta\left(\rho\right)}.
\]
In particular,
\begin{equation}\label{eq:zero-weight-moments}
 \mathcal M V\left(1\right)=\mathcal M V\left(\frac{1}{2}\right)=0,
 \qquad \mathcal M V\left(\rho\right)=1.
\end{equation}
The first vanishing moment ensures smoothness at the origin; the second
cancels the pole of $\zeta\left(2s\right)$.

For $u>0$, put $g\left(u\right)=\sum_{m\geqslant1}V\left(mu\right)$,
and set $g\left(u\right)=0$ for $u\leqslant0$. The sum is locally finite
on $\left(0,\infty\right)$ and vanishes for $u\geqslant\frac{1}{2}$.
Since $\int_0^\infty V\left(u\right)\dd u=0$, the Poisson summation
argument in Lemma~\ref{lem:weight} also applies to this complex-valued $V$.
It gives, for every integer $j\geqslant0$ and every $B>0$,
\begin{equation}\label{eq:zero-weight-flatness}
 g^{\left(j\right)}\left(u\right)\ll_{j,B,\rho}u^B
 \qquad\left(0<u\leqslant1\right).
\end{equation}
Thus $g\in C_c^\infty\left(\mathbb R\right)$, with all derivatives
vanishing at zero. For $\Re s>1$, absolute convergence gives
$\mathcal M g\left(s\right)=\zeta\left(s\right)\mathcal M V\left(s\right)$.
Both sides are entire by \eqref{eq:zero-weight-flatness} and
\eqref{eq:zero-weight-moments}, respectively. Hence the identity holds
for every $s\in\mathbb C$; in particular,
$\mathcal M g\left(\rho\right)=0$.

Define
\[
 w_\rho\left(u\right)
 =u^{-\rho}\int_u^\infty g\left(t\right)t^{\rho-1}\dd t
 \quad\left(u>0\right),
 \qquad w_\rho\left(u\right)=0\quad\left(u\leqslant0\right).
\]
Since $\mathcal M g\left(\rho\right)=0$, a change of variables gives
\[
 w_\rho\left(u\right)
 =-u^{-\rho}\int_0^u g\left(t\right)t^{\rho-1}\dd t
 =-\int_0^1 g\left(uv\right)v^{\rho-1}\dd v.
\]
By \eqref{eq:zero-weight-flatness}, we may differentiate under the integral
for $u>0$. For every integer $j\geqslant0$ and every $B>0$, this gives
\[
 w_\rho^{\left(j\right)}\left(u\right)
 =-\int_0^1 g^{\left(j\right)}\left(uv\right)v^{\rho+j-1}\dd v
 \ll_{j,B,\rho}u^B\int_0^1v^{B+\beta+j-1}\dd v
 \ll_{j,B,\rho}u^B
 \qquad\left(0<u\leqslant1\right).
\]
Also $w_\rho\left(u\right)=0$ for $u\geqslant\frac{1}{2}$.
Thus its zero extension belongs to $C_c^\infty\left(\mathbb R\right)$
with the required support.

Taking Mellin transforms and interchanging the integrals for $\Re s>1$
gives
\begin{equation}\label{eq:zero-weight-mellin}
 \mathcal M w_\rho\left(s\right)
 =\frac{\mathcal M g\left(s\right)}{s-\rho}
 =\frac{\zeta\left(s\right)\mathcal M V\left(s\right)}{s-\rho}.
\end{equation}
Both sides extend to entire functions. At $s=\rho$, the zero of
$\zeta$ cancels the denominator, regardless of its multiplicity.

For $\Re s>1$, the Euler product gives
$\sum_{n\geqslant1}\frac{\lambda\left(n\right)}{n^s}
=\frac{\zeta\left(2s\right)}{\zeta\left(s\right)}$.
Mellin inversion on $\Re s=2$ and \eqref{eq:zero-weight-mellin}
therefore give
\[
 \sum_{n\geqslant1}\lambda\left(n\right)w_\rho\left(\frac{n}{Y}\right)
 =\frac{1}{2\pi\mathrm{i}}\int_{\Re s=2}
   \frac{\zeta\left(2s\right)\mathcal M V\left(s\right)}{s-\rho}
   Y^s\dd s.
\]
The pole at $s=\frac{1}{2}$ is cancelled by
$\mathcal M V\left(\frac{1}{2}\right)=0$. Fix $A>0$ and move the
line of integration to $\Re s=-A$. The only pole crossed is the simple
pole at $s=\rho$, whose residue is
$\zeta\left(2\rho\right)Y^\rho$ by
\eqref{eq:zero-weight-moments}. The rapid decay of $\mathcal M V$
on closed vertical strips and the polynomial growth of $\zeta$ there
justify the shift and bound the integral on the new line by
$O_{A,\rho}\left(Y^{-A}\right)$ for $Y\geqslant1$.
This proves \eqref{eq:zero-weighted-sum}.

Finally, $\zeta\left(2\rho\right)\ne0$ by the Euler product if
$\beta>\frac{1}{2}$, and by the zero-free line $\Re s=1$ if
$\beta=\frac{1}{2}$, since $\gamma\ne0$.
\end{proof}


\begin{thebibliography}{99}

\bibitem{BalogPerelli1998}
A. Balog and A. Perelli,
\emph{On the $L^1$ mean of the exponential sum formed with the M\"obius function},
J. London Math. Soc. (2) \textbf{57} (1998), no.~2, 275--288.

\bibitem{BalogRuzsa2001}
A. Balog and I. Z. Ruzsa,
\emph{On the exponential sum over $r$-free integers},
Acta Math. Hungar. \textbf{90} (2001), no.~3, 219--230.

\bibitem{BrudernEtAl1998}
J. Br\"udern, A. Granville, A. Perelli, R. C. Vaughan and T. D. Wooley,
\emph{On the exponential sum over $k$-free numbers},
Philos. Trans. R. Soc. Lond. Ser. A \textbf{356} (1998), no.~1738, 739--761.

\bibitem{GolinskiiLubinskyNevai2001}
L. Golinskii, D. S. Lubinsky and P. Nevai,
\emph{Large sieve estimates on arcs of a circle},
J. Number Theory \textbf{91} (2001), no.~2, 206--229.

\bibitem{Konyagin1982}
S. V. Konyagin,
\emph{On a problem of Littlewood},
Math. USSR-Izv. \textbf{18} (1982), no.~2, 205--225.

\bibitem{McGeheePignoSmith1981}
O. C. McGehee, L. Pigno and B. Smith,
\emph{Hardy's inequality and the $L^1$-norm of exponential sums},
Ann. of Math. (2) \textbf{113} (1981), no.~3, 613--618.

\bibitem{MontgomeryVaughan2007}
H. L. Montgomery and R. C. Vaughan,
\emph{Multiplicative Number Theory I: Classical Theory},
Cambridge Stud. Adv. Math. \textbf{97}, Cambridge University Press, Cambridge, 2007.

\bibitem{PandeyRadziwill2026}
M. Pandey and M. Radziwi\l\l,
\emph{$L^1$ mean of exponential sums with multiplicative coefficients. I.},
Compos. Math. \textbf{162} (2026), no.~6, 1247--1255.

\bibitem{Pintz1982}
J. Pintz,
\emph{Oscillatory properties of $M(x)=\sum_{n\leqslant x}\mu(n)$. I},
Acta Arith. \textbf{42} (1982), no.~1, 49--55.

\bibitem{Pintz2022}
J. Pintz,
\emph{On the mean value of arithmetic error terms},
Math. Pannon. \textbf{28} (N.S.~2) (2022), no.~1, 58--64.

\bibitem{Yakubovich2015}
S. Yakubovich,
\emph{New summation and transformation formulas of the Poisson, M\"untz,
M\"obius and Voronoi type},
Integral Transforms Spec. Funct. \textbf{26} (2015), no.~10, 768--795.

\end{thebibliography}
\end{document}